\documentclass[12pt]{article}
\usepackage{amsthm,amsmath,amssymb}
\theoremstyle{plain}
\newtheorem{Proposition*}{Proposition}
\newtheorem*{Remark*}{Remark}
\newcommand{\abs}[1]{\lvert#1\rvert}
\newcommand{\FFF}{\mathbb F}            %%% finite field %%%%%%
\newcommand{\NNN}{\mathbb N}
\begin{document}
\title{Permutation groups of prime degree, a quick proof of Burnside's
  theorem} \author{Peter M\"uller} \maketitle

\begin{abstract}
  A transitive permutation group of prime degree is doubly transitive or
  solvable. We give a direct proof of this theorem by Burnside which uses
  neither S-ring type arguments, nor representation theory.
%\footnote{Primary: 20G40, Secondary: none}
\end{abstract}

In this note $\FFF_p$ is the field with $p$ elements, with $p$ a prime. The
following proposition proves Burnside's theorem in a few lines.

\begin{Proposition*}
  Let $U$ be a non-empty, proper subset of $\FFF_p\setminus\{0\}$. Let $\pi$
  be a permutation of $\FFF_p$ such that $i-j\in U$ for $i,j\in\FFF_p$ implies
  $\pi(i)-\pi(j)\in U$. Then there are $a,b\in\FFF_p$ such that $\pi(i)=ai+b$
  for all $i\in\FFF_p$.
\end{Proposition*}
In \cite{Schur:1908} Schur gives a proof of this proposition in two steps.
First he uses a precursor of his S-ring technique to show that if $1\in U$,
then $U$ is a subgroup of $\FFF_p\setminus\{0\}$. In the second step he shows
that $\pi$ is linear. In this note we show that a small modification of his
second step makes the first step unnecessary. See the remarks at the end for
further comments.
\paragraph{Proof of Burnside's theorem.} Let $G$ be a transitive permutation
group on $p$ elements. As $p$ divides the order of $G$, there is an element
$\tau\in G$ of order $p$. Assume that $G$ acts on $\FFF_p$, with $\tau(i)=i+1$
for all $i\in\FFF_p$. Suppose that $G$ is not doubly transitive.  So $G$ has
at least two orbits on the pairs $(i,j)$ with $i\ne j$. On the other hand,
$\tau$ permutes cyclically the pairs $(i,j)$ with constant difference, so
there is a non-empty proper subset $U$ of $\FFF_p\setminus\{0\}$ such that
$\pi(i)-\pi(j)\in U$ for all $\pi\in G$ and $i,j$ with $i-j\in U$. By the
proposition, $G$ is a subgroup of the group of permutations $i\mapsto ai+b$
with $a\in\FFF_p\setminus\{0\},b\in\FFF_p$. In particular, $G$ is solvable.

\paragraph{Proof of the proposition.} By an iterated application of $\pi$ we
see that $i-j\in U$ if and only if $\pi(i)-\pi(j)\in U$. In particular,
replacing $U$ by its complement in $\FFF_p\setminus\{0\}$ preserves the
assumption.  Therefore we may and do assume $\abs{U}\le\frac{p-1}{2}$.

Fix $i\in\FFF_p$. For $u\in U$ we have $(i+u)-i\in U$, hence
$\pi(i+u)-\pi(i)\in U$. As $\pi$ is a permutation, the elements
$\pi(i+u)-\pi(i)$ are different for different $u$. Thus
$\{\pi(i+u)-\pi(i)|\,u\in U\}=U$, hence $\{\pi(i+u)|\,u\in
U\}=\{\pi(i)+u|\,u\in U\}$. In particular, for $w\in\NNN$ we obtain
\[
\sum_{u\in U}\pi(i+u)^w=\sum_{u\in U}(\pi(i)+u)^w.
\]
Let $f(X)\in\FFF_p[X]$ be the polynomial of degree $n\le p-1$ with
$f(i)=\pi(i)$ for all $i\in\FFF_p$. Suppose $wn\le p-1$. Then $\sum_{u\in
  U}f(X+u)^w-\sum_{u\in U}(f(X)+u)^w$ is a polynomial of degree $<p$ which
vanishes identically on $\FFF_p$, thus
\[
\sum_{u\in U}f(X+u)^w-\sum_{u\in U}(f(X)+u)^w=0.
\]
Setting $S(k)=\sum_{u\in U}u^k$, we obtain
\[
\sum_{u\in U}(f(X+u)^w-f(X)^w)=\sum_{k\ge1}\binom{w}{k}S(k)f(X)^{w-k}.
\]
Note that $f(X)^w$ is a polynomial of degree $nw$, so $X^{nw}$ is an
$\FFF_p$-linear combination of the derivatives of $f(X)^w$. Thus we obtain
\[
\sum_{u\in U}((X+u)^{nw}-X^{nw})=\sum_{k\ge1}S(k)g_{w-k}(X),
\]
where $g_\ell(X)$ is a polynomial of degree at most $\ell n$.

Let $r\ge1$ be minimal with $S(r)\ne0$. Then the degree of the right handside
is at most $n(w-r)$. 

Suppose that $r\le nw$. Then the coefficient of $X^{nw-r}$ on the left
handside is (up to a nonzero factor) $S(r)$. Since $S(r)\ne0$, we must have
$nw-r\le n(w-r)$, so $n=1$, and we are done.

It remains to consider $r-1\ge nw$. Suppose we have chosen $w$ maximal with
$nw\le p-1$. Then $p-1<n(w+1)\le 2nw\le 2(r-1)$, so $r>(p+1)/2$.

This shows $S(k)=0$ for $k=1,2,\dots,(p-1)/2$. Therefore $\abs{U}\ge (p+1)/2$
(for instance because the van der Monde matrix of $U$ is not singular; or
because the first $(p-1)/2$ elementary symmetric functions of $U$ vanish, so a
polynomial with zero set $U$ has degree $\ge(p+1)/2$). This contradicts $\abs{U}\le(p-1)/2$.
\begin{Remark*}
  Our proof is similar to the final step in Schur's proof in
  \cite{Schur:1908}. However, the main part of his proof consists in showing
  that if $1\in U$, then $U$ is a subgroup of $\FFF_p\setminus\{0\}$. Thus if
  $1\le k\le w<\abs{U}$, then $S(k)=0$, so $\sum_{u\in
    U}f(X+u)^w=\abs{U}f(X)^w$, which produces a contradiction similarly as
  above.  See also \cite[3.5]{DixMort} for a modern version of this proof.
  
  In \cite{Dress:Burnside} the authors give an S-ring argument to show that
  $U$ is a group. From there they however proceed with geometric arguments,
  and use facts about lacunary polynomials to conclude that $\pi$ is a linear
  function.
  
  Burnside's original proof uses complex character theory, see
  \cite{Burnside:1911}.
  
  The certainly most elegant proof is due to Wielandt, who studies the ring of
  $G$-invariant functions from $\FFF_p$ to $\FFF_p$. See \cite[pages
  273--296]{Wielandt:Coll}, \cite[XII, \S10]{Huppert3}. A very concise and
  streamlined version of Wielandt's proof is contained in \cite[6.7]{LiMuTu}.
\end{Remark*}

\noindent{\sc IWR, Universit\"at Heidelberg, Im Neuenheimer Feld 368,\\
69120 Heidelberg, Germany}\par
\noindent{\sl E-mail: }{\tt Peter.Mueller@iwr.uni-heidelberg.de}
\end{document}